\documentclass[10pt]{amsart}
\usepackage{amsthm,amsfonts,amsmath,amssymb,mathrsfs,enumitem, eucal}
\usepackage[usenames]{color}
\usepackage[colorlinks=true, linkcolor=webgreen, citecolor=webgreen, urlcolor=webbrown]{hyperref}

\definecolor{webgreen}{rgb}{0,.5,0}
\definecolor{webbrown}{rgb}{.5,0,0}

\newtheorem{Thm}{Theorem}[section]

\newtheorem{Lm}[Thm]{Lemma}
\newtheorem{Prop}[Thm]{Proposition}
\newtheorem{Cor}[Thm]{Corollary}
\newtheorem*{state}{Main Theorem}

\theoremstyle{definition}
\newtheorem*{ack}{Acknowledgements}

\theoremstyle{remark}
\newtheorem{introRem}{Remark}
\newtheorem{Rem}[Thm]{Remark}

\numberwithin{equation}{section}

\def\<{\langle}
\def\>{\rangle}

\DeclareMathOperator{\rank}{rank}
\DeclareMathOperator{\trdeg}{trdeg}
\DeclareMathOperator{\gr}{gr}
\DeclareMathOperator{\id}{id}

\def\D{\Delta}
\def\O{\Omega}

\def\M{\mathfrak{m}}
\def\a{\mathfrak{a}}

\def\A{\mathcal{A}}
\def\P{\mathcal{P}}
\def\e{\varepsilon}

\begin{document}
\title[]
{A Poincar\' e-Birkhoff-Witt theorem for Hopf algebras with central Hopf algebra coradical}
\author{Bogdan Ion}
\thanks{Department of Mathematics, University of Pittsburgh, Pittsburgh, PA15260 and University of Bucharest, Faculty of Mathematics and Computer Science, Algebra and Number Theory Research Center, 14 Academiei St., Bucharest, Romania}
\thanks{E-mail address: \nolinkurl{bion@pitt.edu}}
\date{November 8, 2009}%
\subjclass[2000]{14L15, 16W30, 16W70}


\begin{abstract}
We show that over  fields of characteristic zero a Hopf algebra with central Hopf algebra coradical has a PBW basis.
\end{abstract}
\maketitle

\section*{Introduction} 

The Poincar\' e-Birkhoff-Witt theorem is a classical theorem in Lie theory that specifies a special basis for the enveloping algebra of a  Lie algebra. To my knowledge, extensions of this theorem to other classes of algebras depend upon more or less explicit knowledge of the nature of the generators and relations describing the algebras in question, with the notable exception of \cite[Theorem 2]{lebruyn}. 

Let $k$ be a field,  $A$ be a $k$-algebra,  $B$ be a commutative $k$-algebra, and  $B\to A$ a $B$-algebra structure on $A$ (i.e. the image of $B$ inside $A$ is central). For a well-ordered set $\Lambda$ and $\A_\Lambda=\{a_\lambda\}_{\lambda\in \Lambda}$ a subset of $A$ we call 
\begin{equation*}
\P_\Lambda:=\left\{a^{k_1}_{\lambda_1}\dots a^{k_n}_{\lambda_n}~\Bigl\lvert~ n, k_1,\dots, k_n\geq 0,~ \lambda_1<\dots<\lambda_n\right\}
\end{equation*}
the PBW set associated to the generating set $\A_\Lambda$. We say that $A$ has a PBW basis over $B$ if there exist a well-ordered set $\Lambda$ and a generating set $\A_\Lambda$ such that the PBW set associated to $\A_\Lambda$ is a basis, called PBW basis, of $A$ as a $B$-module.

We will be concerned with Hopf $k$-algebras $H$ with central Hopf algebra coradical $H_0$. Note that if such an algebra is finite dimensional then it is, according to Corollary \ref{finitefiltered},  automatically equal to its coradical. For Hopf algebras with proper coradical, hence infinite dimensional, PBW bases start playing an important role. Our main result is the following.
\begin{state}\label{state}
Let $k$ be a field of characteristic zero and let $H$ be a Hopf $k$-algebra such that its coradical $H_0$ is a central  Hopf subalgebra. Then $H$ has a PBW basis over $H_0$.
\end{state}
The strategy for the proof, as in \cite[Theorem 2]{lebruyn}, is to show that $\gr H$, the associated graded algebra of $H$ with respect to the coradical filtration, is a polynomial algebra with coefficients in $H_0$. For an extension of these results to symmetrically braided Hopf algebras we refer to \cite{symPBW}.

Before proceeding with the proof of these results let us comment on the nature of the constraints present in the hypothesis as well as its relationship with other PBW theorems in the literature. 

\begin{introRem} The notion of PBW basis in this paper is stronger than the one generally used in the literature. The restricted PBW bases are those for which the exponents $k_1,\dots, k_n$ from the definition are non-negative integers that satisfy explicitly specified upper bounds (possibly infinity). As such, it is possible for finite dimensional algebras to have restricted PBW bases and this is indeed the case for many important examples such as enveloping algebras of restricted Lie algebras and Lusztig's small quantum groups \cite{lusztig} and their deformations (see, e.g., \cite{andrus}).
\end{introRem}
\begin{introRem} If one allows restricted PBW bases then the condition that $H_0$ is a central Hopf subalgebra of $H$ can be weakened, at least for some classes of Hopf algebras. The most general results of this type currently in the literature are due to Kharchenko \cite{khar1, khar2} and they cover  the class of the so-called character Hopf algebras: Hopf algebras generated by an abelian subgroup of group-like elements and by a finite set of skew-primitive elements such that the adjoint action of this abelian group on the skew-primitive generators is given by multiplication with a character. It seems that most if not all finite dimensional pointed Hopf algebras with abelian coradical fall within this class \cite{andrus}.
\end{introRem}
\begin{introRem} Other algebras for which PBW theorems are of interest are Hopf algebras in braided categories. Again, there are numerous examples of pointed irreducible Hopf algebras (therefore with central coradical) in categories of Yetter--Drinfel'd modules that are finite dimensional and therefore do not admit a PBW basis in our sense but have a restricted PBW basis. One possible reformulation of the result in \cite{khar1} is that braided tensor algebra quotients of the tensor algebra of a Yetter--Drinfel'd module (over an abelian group) with diagonal braiding admit restricted PBW bases. This result was extended to cover also triangular braidings \cite{ufer}.
\end{introRem}
\begin{introRem}
The combinatorial arguments in \cite{khar1, khar2, ufer} require fairly explicit knowledge of the generators and relations describing the algebras in question especially when it comes to the description of the generating set of the PBW basis. On the other hand, as it will become clear from the proof, the core ingredient underpinning our approach is  smoothness, and the generating set of our PBW set acquires a fairly simple conceptual characterization.
\end{introRem}
\begin{introRem}
The ideas used here can be also employed to prove the corresponding result over a field $k$ of characteristic $p>0$. The necessary modifications of the statements are the following: the exponents $k_1,\dots, k_n$ from the definition of $\P_\Lambda$ have to be strictly less than $p$;  the algebra $\gr H$ is the quotient of a polynomial algebra with coefficients in $H_0$ by the ideal generated by the $p^{th}$ powers of the variables. The only additional technical observation is the fact that the affine group scheme corresponding to the coinvariants $R$ (see Section \ref{coinv}) is of height one \cite[$\S$ 11.4]{waterhouse}.
\end{introRem}
\begin{ack}
 Work supported in part by NSF grant DMS--0536962 and CNCSIS grant nr. 24/28.09.07 (Groups, quantum groups, corings, and representation theory).  Some of the results of this paper are  to some extent contained in \cite{ion}. I thank \c Serban Raianu for suggesting the topic of \cite{ion} and for bringing \cite[Theorem 2]{lebruyn} to my attention.  
\end{ack}

\section{Preliminaries}\label{prelim}
\subsection{} Throughout the paper $k$ will denote a field of characteristic zero and $\bar{k}$ will denote its algebraic closure. Unless otherwise noted all tensor products are over $k$.
When referring to a Hopf algebra, we use $\D$, $\e$, and $S$ to denote its comultiplication, counit, and antipode, respectively. As general notation we will use capital letters to denote Hopf algebras and will use  the corresponding small case gothic letters to denote their augmentation ideal; e.g. the augmentation ideal of the Hopf algebra $A$ is denoted by $\a$.  All the results in Section \ref{prelim} are well-known. 

\subsection{}\label{polmorph} Assume $A$ is a $k$-algebra, $B\subseteq A$ is a $k$-subalgebra, and $S\subseteq A$ is a subset. If $A$ is commutative then denote by 
\begin{equation}
p_{B, S}: B[X_s]_{s\in S}\to A
\end{equation}
the ring morphism acting as identity on $B$ and sending each $X_s$ to the element $s$. We will drop the reference to $B$ and $S$ if it is unambiguous from the context.

\subsection{} Henceforward $G$ will denote an affine algebraic group scheme over $k$ and $A$ will denote the finitely generated commutative Hopf $k$-algebra representing it. 
Let $\pi_\circ A$ denote the largest separable $k$-subalgebra of $A$. It is a finite dimensional Hopf subalgebra of $A$ that represents the \' etale group scheme $\pi_\circ G$, the component group of $G$. 

The kernel $G^\circ$ of the canonical map $G\to\pi_\circ G$ is the connected component of the identity in $G$. It is itself an affine algebraic groups scheme; we denote by $A^\circ$ the finitely generated commutative Hopf $k$-algebra representing it (a Hopf algebra quotient of $A$). The construction of $\pi_\circ G$ and  $G^\circ$ commutes with base change.

\subsection{}\label{components} 
The algebra $\pi_\circ A$ is a product of finite field extensions of $k$ and $\pi_\circ A\otimes \bar{k}$ is isomorphic to $\bar{k}\times\cdots\times\bar{k}$ where the number of factors equals the number of connected components of $G(\bar{k})$. It is clear from this description that the $\bar{k}$-algebra $\pi_\circ A\otimes \bar{k}$ is spanned by idempotents of $A\otimes\bar{k}$.

Using the factorization of $\pi_\circ A$ into fields we can write $A=\oplus_{i=0}^s f_i A$, with $f_i$ idempotents in $A$. The counit $\e$ must vanish on all but one of these idempotents; we may assume that $f_0$ is the idempotent on which $\e$ takes value 1. With this notation, the canonical map $G^\circ\to G$ corresponds to 
\begin{equation}
\pi_A: A\to A/(\a\cap\pi_\circ A)A=A/(1-f_0)A = A^\circ\simeq f_0A 
\end{equation}
 The algebras $f_iA$ represent the connected connected components of $G$. They are generally not isomorphic as algebras unless $k=\bar{k}$. 
\subsection{} Let $\O_A$ be the universal module of differentials of $A$ and let $d_A: A\to \O_A$ the universal derivation. The pair $(\O_A, d_A)$ is unique up to isomorphism. Among other things, the construction of $\O_A$ commutes with base change and localization.

Let  $\pi: A=k\cdot 1\oplus \a\to \a\to \a/\a^2$. With this notation,
\begin{equation}
(\O_A, d_A)\simeq \left(A\otimes \a/\a^2, ({\rm id}_A\otimes \pi)\circ\D\right)
\end{equation} 
In particular, $\O_A$ is a free $A$-module of rank equal to the $k$-dimension of $\a/\a^2$. 
\subsection{} Since we work over fields of characteristic zero the following results are available. For a proof of the first theorem, generally attributed to Cartier, see \cite[$\S$ 11.4]{waterhouse}, and  for a proof of the second theorem see \cite[$\S$ 11.6]{waterhouse}.
\begin{Thm}\label{cartier}
The Hopf algebra $A$ is reduced.
\end{Thm}
One important consequence is the following.
\begin{Cor}\label{finite=etale}
All finite group schemes over $k$ are \' etale.
\end{Cor}
\begin{Thm}\label{smooth}
The affine algebraic group scheme $G$ is smooth.
\end{Thm}
If apply Theorem \ref{cartier} for $G^\circ$, which is connected, we obtain that $A^\circ$ is an integral domain; we denote by $A^\circ_{(0)}$ its field of fractions. Theorem \ref{smooth} now reads
\begin{equation}\label{omegarank}
\rank_A\O_A=\dim_k\pi_\circ A\cdot\trdeg_kA^\circ_{(0)}
\end{equation}
We can actually say more.
\begin{Thm}\label{omegagenerators} 
Let $x_1,\dots, x_n$ be elements of $\a$ whose classes modulo $\a^2$ form a $k$-basis of $\a/\a^2$. Then, $d_{A^\circ_{(0)}}(x_1), \dots, d_{A^\circ_{(0)}}(x_n)$ form a basis for $\O_{A^\circ_{(0)}}$ as a ${A^\circ_{(0)}}$-vector space. 
\end{Thm}
The proof is contained in \cite[$\S$ 13.5]{waterhouse}. Using translation in $G(k)$, which by Theorem \ref{smooth} is smooth, the same is true if we replace the identity element (i.e. the augmentation ideal) with any other point of $G(k)$  (i.e. the kernel of an algebra morphism $A\to k$). Another useful result is the following
\begin{Thm}\label{fgen}
Let $B$ be a Hopf subalgebra of $A$. Then $B$ is finitely generated.
\end{Thm}
For the proof we refer to \cite[$\S$ 14.3]{waterhouse}.
\subsection{}  Let $k\subseteq L$ be a finitely generated field extension. Regard $L$  as a finitely generated commutative $k$-algebra; we use the notation $L_k$ to stress this.
Let $\O_{L_k}$, and  $d_{L_k}: L\to \O_{L_k}$, be the universal module of differentials of $L_k$, and the universal derivation, respectively. 

Let $x_1,\dots, x_n$ be elements of $L$ such that $d_{L_k}(x_1),\dots, d_{L_k}(x_n)$ is an $L$-basis of $\O_{L_k}$. 
\begin{Thm}\label{algind}
With the notation above, $x_1,\dots, x_n$ are algebraically independent over $k$ and $k(x_1,\dots, x_n)\subseteq L$ is a finite field extension.
\end{Thm}
For the proof see \cite[$\S$ 11.5]{waterhouse}.
\subsection{}\label{induction} The following result is again well-known; we refer to \cite[$\S$ 3.3]{waterhouse} for a proof.
\begin{Thm}\label{inductive}
Any commutative Hopf $k$-algebra is a directed union of finitely generated Hopf subalgebras.
\end{Thm}
\subsection{} Consider now $H$ an arbitrary Hopf $k$-algebra and let $0=:H_{-1}\subseteq H_0\subseteq H_1\subseteq \dots \subseteq H$ its coradical filtration. Denote by 
\begin{equation}
\gr H=\oplus_{n\geq 0}(\gr H)(n)=\oplus_{n\geq 0} H_n/H_{n-1}
\end{equation}
the associated graded coalgebra. If $B\subseteq H$ is a Hopf subalgebra then the coradical filtration of $B$ is the one induced from the coradical filtration of $H$.

 If the coradical $H_0$ is a Hopf subalgebra of $H$ then the coradical filtration is a Hopf algebra filtration and $\gr H$ is a graded Hopf algebra. Furthermore, in such a case, the coradical filtration of $\gr H$ is the filtration given by degree and $\gr(\gr H)=\gr H$.

The the proof of the following result we refer to \cite[Proposition 1.2]{montgomery}.
\begin{Prop}\label{grcomm}
Let $H$ be a Hopf algebra such the coradical $H_0$ is a central Hopf subalgebra. Then $\gr H$ is a commutative Hopf algebra. 
\end{Prop}

\section{The main result} 

\subsection{} We start with a few observations on commutative Hopf $k$-algebras endowed with a Hopf algebra grading. 
\begin{Lm}\label{lemma}
Let  $A=\oplus_{i\geq 0}A(i)$ be a finitely generated commutative Hopf $k$-algebra endowed with a Hopf algebra grading. Then $\pi_\circ A=\pi_\circ A(0)$. 
\end{Lm}
\begin{proof} First remark that by Theorem \ref{fgen} the Hopf algebra $A(0)$ is also finitely generated. Any separable subalgebra of $A(0)$ is necessarily a separable subalgebra of $\pi_\circ A$, so $\pi_\circ A(0)\subseteq \pi_\circ A $. To settle the claim it would be enough to show that they have the same dimension. Since the construction of both subalgebras commutes with base change it is safe to assume that $k=\bar{k}$.

In such a case, $\pi_\circ A$ is spanned by idempotents and it is therefore enough to show that all idempotents in $A$ are necessarily in $A(0)$. Indeed, let
\begin{equation}
e=e_{i_1}+\cdots+e_{i_k}
\end{equation}
be an idempotent and $e_{i_1}, \dots, e_{i_k}$, $i_1<\dots<i_k$ be its non-zero homogeneous  components. Then, $e^2_{i_k}$ is the homogeneous  component of degree $2i_k$ of $e^2=e$. From Theorem \ref{cartier}  we know that $A$ is reduced. Hence we must have $i_k=0$. Therefore, $e$ is an element of $A(0)$.
\end{proof}
\begin{Cor}\label{finitegraded}
Let  $A=\oplus_{i\geq 0}A(i)$ be a finite dimensional commutative Hopf $k$-algebra endowed with a Hopf algebra grading. Then $A=A(0)$.
\end{Cor}
\begin{proof} From Corollary \ref{finite=etale} we know that $A=\pi_\circ A$ and $A(0)=\pi_\circ A(0)$. The desired statement follows then from the above Lemma.
\end{proof}
\begin{Cor}\label{finitefiltered}
Let $H$ be a finite dimensional Hopf $k$-algebra with central Hopf algebra coradical. Then $H$ is cosemisimple.
\end{Cor}
\begin{proof}
According to Proposition \ref{grcomm} the ring $A:=\gr H$ is a commutative Hopf $k$-algebra endowed with a Hopf algebra grading such that $A(0)$ is the coradical of $H$. But then $A=A(0)$ which implies our claim.
\end{proof}
\subsection{}\label{coinv}
Again, let  $A=\oplus_{i\geq 0}A(i)$ be a finitely generated commutative Hopf $k$-algebra endowed with a Hopf algebra grading. The projection $\varpi: A\to A(0)$ is a Hopf algebra retraction of the inclusion $\iota: A(0)\to A$. In such a situation one can define the subalgebra of coinvariants of $\varpi$
\begin{equation}\label{R1}
R:=\{a\in A~|~(\id\otimes \varpi)\circ\Delta(a)=a\otimes 1 \}=\left\{\sum a_1S(\varpi(a_2))~|~a\in A\right\}
\end{equation}
where above we used the Sweedler notation $\Delta(a)=\sum a_1\otimes a_2$. Let $R(i):=R\cap A(i)$. It is easy to verify (see, e.g. \cite[Lemma 2.1]{andrus2}) that $R=\oplus_{i\geq 0}R(i)$ is a graded algebra and that the canonical map $A(0)\otimes R\to A$, $a\otimes r\mapsto ar$ which is an algebra isomorphism by \cite[Theorem 3]{radford} is a compatible with the grading.

Denote $A_+:=\oplus_{n\geq 1}A(n)$ and $R_+:=\oplus_{n\geq 1}R(n)$. The $A(0)$-module isomorphism $A(0)\otimes R_+\simeq  A_+$ descends to a graded vector space isomorphism
\begin{equation}
R_+\simeq \frac{A_+}{\a(0)A_+}
\end{equation}
The following result plays a crucial role in the proof of Proposition \ref{prop1}.
\begin{Lm}\label{lemmaR}
With the notation above, there exists a canonical vector space isomorphism
$$
\frac{R_+}{R_+^2}\simeq \frac{A_+}{\a(0)A_++A_+^2}
$$ 
\end{Lm}
\begin{proof}
Consider the canonical $k$-linear map $$\gamma: R_+\simeq \frac{A_+}{\a(0)A_+}\to \frac{A_+/\a(0)A_+}{(\a(0)A_++A_+^2)/\a(0)A_+}\simeq \frac{A_+}{\a(0)A_++A_+^2}$$
We argue that the kernel of $R^2_+=\ker\gamma$. The direct inclusion is obvious. For the reverse inclusion, let $r$ be in $\ker\gamma$. Since $\a(0)A_++A_+^2=\a(0)R_++A(0)R_+^2$ there exist elements $a_i\in A(0)$, $r_j, s_k\in R_+$ such that 
$$
r-\sum a_ir_js_k\in \a(0)A_+
$$
Keeping in mind that $a_i-\e(a_i)\in\a(0)$ we deduce that 
$$
r-\sum \e(a_i)r_js_k\in \a(0)A_+
$$
from which our claim immediately follows.
\end{proof}
Let $A=\oplus_{i=0}^s f_i A$ be the decomposition discussed in Section \ref{components}. The elements $f_i$ are idempotents so, according to Lemma \ref{lemma}, they are elements of $A(0)$. Recall that  $\e(f_0)=1$ and $\e(f_i)=0$ for $i\neq 0$.
\begin{Lm}\label{lemmaR2}
With the notation above, $\ker \pi_A\cap R=0$.
\end{Lm}
\begin{proof}
We will show that $f_iA\cap R=0$ for any $i\neq 0$. Let $a\in f_iA\cap R$. Then,  $a=f_ia$  and  $a=\sum a_1S(\varpi(a_2))$. Therefore, 
\begin{align*}
a&=\sum f_{i,1}a_1S(\varpi(a_2))S(\varpi(f_{i,2}))\\
&= \e(f_i)\sum a_1S(\varpi(a_2))\\
&=0
\end{align*}
which is our claim.
\end{proof}
\subsection{} We now prove a first version of our main result. 
\begin{Prop}\label{prop1}
Let $A=\oplus_{i\geq 0}A(i)$ be a  finitely generated commutative Hopf $k$-algebra  endowed with a Hopf algebra grading.  Assume that $A$ is an integral domain.  Then there exist a finite set $\A_\Lambda$ of $A(0)$-coinvariant homogeneous elements of positive degree  such that  $p_{A(0),\A_\Lambda}$ is an algebra isomorphism.
\end{Prop}
\begin{proof}  The subspace 
\begin{equation}
\M:=\a(0)\oplus A_+
\end{equation}
is a graded maximal ideal of $A$. We can also write $A=k\cdot 1\oplus \M$. Of course, 
\begin{equation}\label{eq1}
\frac{\M}{\M^2}\simeq \frac{\a(0)}{\a(0)^2}\oplus\frac{A_+}{\a(0)A_++A_+^2}
\end{equation}
so we may chose $x_1,\dots, x_m$ elements of $\M$ such that $x_i$, $1\leq i\leq n$ are homogeneous elements of $A_+$ and $x_i$, $n<i\leq m$ are elements of  $\a(0)$ and such that their classes modulo $\M^2$ form a $k$-basis of $\M/\M^2$. In addition, as a consequence of Lemma \ref{lemmaR}, we can arrange that the elements $x_i$, $1\leq i\leq n$ are homogeneous elements of $R_+$ whose classes form a basis of the vector space $R_+/R_+^2$.

Keeping in mind that $R$ is a graded algebra with $R(0)=k\cdot 1$ we deduce that it is generated as a $k$-algebra by $x_i$, $1\leq i\leq n$. Therefore, $A$ is generated as a $A(0)$-algebra by $x_i$, $1\leq i\leq n$. Let us denote by 
\begin{equation}\label{eq2}
p: A(0)[X_1,\dots, X_n]\to A
\end{equation}
the canonical morphism that acts identically on $A(0)$ and send each $X_i$ to $x_i$. We will show that  $p$ is an isomorphism.

The fact that $p$ is surjective implies that 
\begin{equation}\label{eq5}
A_{(0)}=A(0)_{(0)}(x_1,\dots, x_n)
\end{equation}
Theorem \ref{fgen} assures that $A(0)$ is itself finitely generated. Theorem \ref{omegagenerators} applied, on one hand, to $A(0)$ and $\a(0)$, and on the other hand to $A$ and $\M$ implies that $d_{A(0)_{(0)}}(x_i)$, $n< i\leq m$ is a ${A(0)_{(0)}}$-basis of $\O_{A(0)_{(0)}}$, and $d_{A_{(0)}}(x_i)$, $1\leq i\leq m$ is a ${A_{(0)}}$-basis of $\O_{A_{(0)}}$. Applying now Theorem \ref{algind} we obtain that, on one hand,  $x_i$, $n<i\leq m$ are algebraically independent over $k$ and 
\begin{equation}\label{eq3}
K:=k(x_{n+1},\dots, x_m)\subseteq A(0)_{(0)}
\end{equation}
is a finite extension, and on the other hand $x_i$, $1\leq i\leq m$ are algebraically independent over $k$ and 
\begin{equation}\label{eq4}
E:=k(x_{1},\dots, x_m)\subseteq A_{(0)}
\end{equation}
is a finite extension.

Therefore, we have $K\subseteq E$ a purely transcendental extension with transcendence basis $x_i$, $1\leq i\leq n$, and a finite extension $K\subseteq A(0)_{(0)}$.  Also, from \eqref{eq5} the field $A_{(0)}$ is the compositum of $E$ and $A(0)_{(0)}$. It is easily proved (by induction on $n$) that in such a case, the extension $A(0)_{(0)}\subseteq A_{(0)}$ is purely transcendental with transcendence basis $x_i$, $1\leq i\leq n$.

Denote by $\beta$ the canonical inclusion of $A$ into its field of fractions. The algebra morphism
\begin{equation}
A(0)[X_1,\dots, X_n]\xrightarrow[\quad\quad]{p} A \xrightarrow[\quad\quad]{\beta} A_{(0)}=A(0)_{(0)}(x_1,\dots, x_n)
\end{equation}
acts as identity on $A(0)$ and sends each $X_i$ to $x_i$, $1\leq i\leq n$. Since $x_i$, $1\leq i\leq n$ are transcendental over $A(0)_{(0)}$ we obtain that $\beta\circ p$ must be injective and, in consequence $p$ must be injective. But, from \eqref{eq2} we know that $p$ is surjective. Therefore, $p$ is an algebra isomorphism.
\end{proof}
\begin{Rem} In the above proof, if we assign to each variable $X_i$ the degree of the element $x_i$ in $A$ then $p$ becomes an isomorphism of graded rings.  
\end{Rem}
\begin{Rem} In the statement of Proposition \ref{prop1} we can replace the augmentation ideal with any maximal ideal of $A(0)$. 
\end{Rem}
\subsection{}  We now remove the integral domain condition in Proposition \ref{prop1}.
\begin{Thm}\label{thm1}
Let $A=\oplus_{i\geq 0}A(i)$ be a finitely generated commutative Hopf $k$-algebra  endowed with a Hopf algebra grading. Then there exist a finite set $\A_\Lambda$ of $A(0)$-coinvariant homogeneous elements of positive degree  such that  such that  $p_{A(0),\A_\Lambda}$ is an algebra isomorphism.
\end{Thm}
\begin{proof} 
Let $A=\oplus_{i=0}^s f_i A$ be the decomposition discussed in Section \ref{components}. Remark first that $\pi_\circ A=\pi_\circ A(0)$ implies that 
\begin{equation}\label{ker1}
(\a\cap\pi_\circ A)A=(\a(0)\cap\pi_\circ A(0))A
\end{equation}
and furthermore, 
\begin{equation}\label{ker2}
\left((\a\cap\pi_\circ A)A\right)\cap A(0)=(\a(0)\cap\pi_\circ A(0))A(0)
\end{equation}
or, equivalently, $A^\circ(0)$ and $A(0)^\circ$ are canonically isomorphic and which we will tacitly identify henceforth. Denote by $R^\circ$ the coinvariant algebra of $A^\circ$ with respect to $A(0)^\circ$. From \eqref{R1} we obtain that $\pi_A(R)=R^\circ$. Furthermore, from Lemma \ref{lemmaR2} we know that $\pi_A$ is injective on $R$, hence $R$ and $R^\circ$ are canonically identified.

From the proof of Proposition \ref{prop1} we know that $R^\circ$ is a polynomial $k$-algebra. Our claim now follows from the fact that the algebras $A$ and $A(0)\otimes R$ are canonically isomorphic.
\end{proof}
\begin{Rem} In the above proof, if we assign to each variable $X_i$ the degree of the element $x_i$ in $A$ then $p$ becomes an isomorphism of graded rings.  
\end{Rem}
\subsection{}  The hypothesis in Theorem \ref{thm1} that $A$ is finitely generated can be removed under certain conditions. We only include here the statements, all the claims following from the fact that all constructions from the proof of Proposition \ref{prop1} and Theorem \ref{thm1} are compatible with the direct limits. 
\begin{Thm}\label{thm2}
Let $A=\oplus_{i\geq 0}A(i)$ be a commutative Hopf $k$-algebra  endowed with a Hopf algebra grading such that there exist finitely generated graded Hopf algebras $A_\alpha$ such that $A=\varinjlim A_\alpha$ as graded Hopf algebras. Then there exist a set $\A_\Lambda$ of $A(0)$-coinvariant homogeneous elements of positive degree  such that  such that  $p_{A(0),\A_\Lambda}$ is an algebra isomorphism.
\end{Thm}
\begin{Rem} If we assign to each variable $X_\lambda$ the degree of the element $x_\lambda$ in $A$ then $p$ becomes an isomorphism of graded rings.  
\end{Rem}
\subsection{} We are now ready to complete the proof of our main result. 
\begin{proof}[Proof of Main Theorem]  
According to Proposition \ref{grcomm} the ring $A:=\gr H$ is a commutative Hopf $k$-algebra endowed with a Hopf algebra grading whose associated filtration is the coradical filtration of $\gr H$. Let us argue that $A$ satisfies the hypothesis of Theorem \ref{thm2}. By Theorem \ref{inductive} $A$ is a directed union of finitely generated Hopf $k$-algebras $B_\alpha$. Since the coradical filtration on each $B_\alpha$ is the filtration induced by the coradical filtration of $A$ we obtain that $\gr A$ (which is in fact isomorphic to $A$) is the directed union (as graded algebras) of the graded finitely generated Hopf $k$-algebras $A_\alpha:=\gr B_\alpha$, proving our claim.

We can now invoke Theorem \ref{thm2} and its proof to find homogeneous elements of positive degree $\{\bar{x}_\lambda\}_{\lambda\in\Lambda}\subset A$ with the specified properties. Fix now a well-ordering $\leq$ of the set $\Lambda$ and choose representatives $\{{x}_\lambda\}_{\lambda\in \Lambda}=:\A_\Lambda$ in $H$ for the elements $\{\bar{x}_\lambda\}_{\lambda\in \Lambda}$. Denote the PBW set associated to $\A_\Lambda$ by
\begin{equation}
\P_\Lambda=\left\{x^{k_1}_{\lambda_1}\dots x^{k_n}_{\lambda_n}~\Bigl\lvert~ n, k_1,\dots, k_n\geq 0,~ \lambda_1<\dots<\lambda_n\right\}
\end{equation}
All the required properties of $\P_\Lambda$ then follow directly from the corresponding facts in Theorem \ref{thm2}.
\end{proof}


\begin{thebibliography}{10}

\bibitem{andrus2}
{\sc N. Andruskiewitsch and H.-J. Schneider}, Lifting of quantum linear spaces and pointed Hopf algebras of order $p^3$. {\it J. Algebra} {\bf 209} (1998) no. 2, 658--691.

\bibitem{andrus}
{\sc N. Andruskiewitsch and H.-J. Schneider}, On the classification of finite-dimensional pointed Hopf algebras.  {\it Ann. of Math. (2)} {\bf 171} (2010) no. 1, 375--417.

\bibitem{ion}
{\sc B. Ion}, PBW theorems. (Romanian) License diploma thesis, University of Bucharest, 1997.

\bibitem{symPBW}
{\sc B. Ion}, Relative PBW type theorems for symmetrically braided Hopf algebras. {\it arXiv: 1004.5022}.

\bibitem{khar1}
{\sc V. K. Kharchenko}, A quantum analogue of the Poincar\' e-Birkhoff-Witt theorem. {\it Algebra Log. } {\bf 38} (1999) no. 4, 476--507.

\bibitem{khar2}
{\sc V. K. Kharchenko}, PBW-bases of coideal subalgebras and a freeness theorem. {\it Trans. Amer. Math. Soc.} {\bf 360} (2008) no. 10, 5121--5143.

\bibitem{lebruyn}
{\sc L. Le Bruyn}, Lie stacks and their enveloping algebras.  {\it Adv. Math.} {\bf 130} (1997), no. 1, 103--135.

\bibitem{lusztig}
{\sc G. Lusztig}, Finite-dimensional Hopf algebras arising from quantized universal enveloping algebra. {\it J. Amer. Math. Soc.} {\bf 3} (1990), no. 1, 257--296. 

\bibitem{montgomery}
{\sc S. Montgomery}, Some remarks on filtrations of Hopf algebras. {\it Comm. Algebra} {\bf 21} (1993), no. 3, 999--1007.

\bibitem{radford}
{\sc D. E. Radford}, The structure of Hopf algebras with a projection. {\it J. Algebra} {\bf 92} (1985), no. 2, 322--347.

\bibitem{ufer}
{\sc S. Ufer}, PBW bases for a class of braided Hopf algebras. {\it J. Algebra} {\bf 280} (2004), no. 1, 84--119.

\bibitem{waterhouse}
{\sc W. C. Waterhouse},  Introduction to affine group schemes. {\it Graduate Texts in Mathematics} {\bf 66}. Springer-Verlag, New York-Berlin, 1979. 

\end{thebibliography}
\end{document}